\documentstyle{amsppt}

\def\ss{\smallskip}

\def\O{\Omega}

\def\e{\epsilon}

\def\l{\lambda}

\def\no{\noindent}
\def\bs{\bigskip}

\def\g{\gamma}

\def\b{\beta}
\def\a{\alpha}
\def\n{\nabla}

% "Sparar" check-accent som /hacek, innan /v omdefinieras.
\nologo
\magnification=1200

\topmatter

\document

\title 
A Higher Order Weierstrass Approximation Theorem
-- a new proof
\endtitle

\author 
Andreas Wannebo
\endauthor

\abstract
The theorem studied is known before. Here is given a new proof.
The proof has been part of a course material on Sobolev space theory
with a special kind of outlook.
The proof here is in accordance to this goal.
At the same time several ideas of interest are shown 
that can be of general use.
\endabstract

\endtopmatter

\heading
Introduction
\endheading

The definition of Sobolev space used here is for reasons to obvious later on
that of the closure in the approprate Sobolev space norm
of the set of functions with $m$ times continuous derivatives 
on some open subset of $N$-dimensional Euclidean space.

The Weierstrass approximation theorem is stated below and concerns 
approximation of contionuous functions by polynomials. The higher
order Weierstrass approximation theorem concerns approximation
of $m$ times continuously differetiable functions by polynomials 
in such a way that
all derivatives up to order $m$ are approximated simultanuously.
The theorem is stated below. The theorem is previuosly known.
Here is given a very interesting proof that makes use of several ideas. 
These and their interplay displayed is the real goal!
That a straightforward proof has less complexity is beside the point here. 
\ss

The theorem holds for a cube as domain and
with the open bounded domain denoted $\O$ in $\text{\bf R}^N$,
then, if the domain allows the Sobolev functions considered 
to be extended to the full space, then the theorem follows for this
domain as well. The argument is simply that take cube enclosing the
domain and the do the extension of the function to full 
space and then restrict
to the cube. Then use the approximation by polynom for the cube as in 
the theorem. Conclude by taking the restriction of the polynomial
to the domain.
\ss

A very good extension theorem is that of P.W. Jones [1]. The domains
allowing extension by his method are called uniform domains or 
also by the name $(\e,\delta)$-domains. 
Domains with self-similar fractal boundary are simple examples of
uniform domains.
There is a point here that the Jones extension theorem
also allows for $p=\infty$ as well.
\ss

A consequence is that polynomials are dense 
in integer order Sobolev spaces with $1\le p\le \infty$ when the 
domain is bounded and is a so called uniform domain. 
\vfill
\eject

\heading
Theorem and Proof
\endheading

{\bf Definition:} Let $\O$ be open in $\text{\bf R}^N$.
For $u$ a ``nice'' real function on $\O$ the Sobolev space norm is
$$
||u||_{W^{m,p}(\O)}=\sum _{|\a|\le m}||D^{\a}u||_{L^p(\O)}
\sim \sum _{|\a|\le m}(||D^{\a}u||^p_{L^p(\O)})^{1\over p}.
$$
\bs

{\bf Definition:} Let $\O$ be open in $\text{\bf R}^N$. 
Define $C^m(\O)$ as the $m$ times continuously differentiable
real-valued functions on $\O$.
\bs

{\bf Notation:} We use $\n^ku$, the $k$-gradient, in integral 
expressions as follows
$$
||\n^ku||_{L^p(\O)}=\sum_{|\a|=k}||D^{\a}u||_{L^p(\O)}.
$$
\bs

{\bf Notation:} Let ${\Cal P}$ be all polynomials in $\text{\bf R}^N$
and ${\Cal P}_k=\{P\in {\Cal P}:\n^{k+1}P=0\}$.
\bs

Then some easy results without proofs.
\bs

{\bf Theorem:} A {\bf Poincar\'e inequality}, dim general, order one
and general $p$.
\ss

Let $1\le p\le \infty$ and  $Q=[0,1]^N$. Let $u\in C^0(Q)$ and 
$D_1u\in C^0(Q)$. Furthermore $u|_{x_1=0}=0$. 
Then
$$
||u||_{L^p(Q)}\le ||D_1u||_{L^p(Q)}.
$$
So, here $Q$ is a unit cube and $D_1$ is the derivative in the $x_1$-direction.
\bs

{\bf Corollary:} The {\bf Standard Poincar\'e Inequality.}
\ss

Let $u\in C^m_0(\O)$ for $\O$ bounded, open subset of $\text{\bf R}^N$
and $|\a|<m$. Let $A=A_{N,m,p,\O}$,
$$
\sum_{|\a|<m} ||D^{\a}u||_{L^p(\O)}\le A||\n^mu||_{L^p(\O)}.
$$
\bs

{\bf Corollary:} A more detailed Poincar\'e inequality.
\ss

Let $Q=[0,1]^N$ and let $\a$ and $\b$ be
multi-indices that are partially ordered by ``$<$'' in the natural way.
Let $|\a|=t$, $|\b|=m$ and $\a<\b$.
Assume that  multi-indices $\{\g_k\}_t^m$ exist
such that $\g_t=\a$ and $\g_m=\b$ and also
$\g_t<\g_{t+1}<\dots<\g_{m-1}<\g_m$.

Let $u\in C^m(\bar Q)$ and 
$$
\big(D^{\g_k}u|_{x_{\g_k-\g_{k-1}}=0}\big)=0
$$
or
$$
\big(D^{\g_k}u|_{x_{\g_k-\g_{k-1}}=1}\big)=0.
$$
Then for $A=A_{N,m,t,p,Q}$,
$$
||D^{\a}u||_{L^p(Q)}\le A||D^{\b}u||_{L^p(Q)}.
$$
\bs

{\bf Theorem:} (Weierstrass Approximation Theorem.)
\bs

Let $\O\subset\text{\bf R}^N$ be connected, open and bounded.
Then $\Cal P$ is dense in $L^{\infty}$-norm in $C^0(\bar \O)$.
\bs

This will be a tool for proving a higher order version.
\bs

The proof here is new.
\bs

{\bf Theorem:} A Higher Order {\bf Weierstrass Approximation theorem.}
\bs

Let $m>0$, integer, and let $Q=[0,1]^N$. Then 
${\Cal P}$ is dense in $W^{m,\infty}$--norm in $C^m(\bar Q)$.
\bs

The difficulty with the use of Weierstrass approximation theorem here
is that there are many partial derivatives. They shall all be approximated by 
the partial derivatives of only one polynomial. 
-- A trick (=method) is called for. 

First an observation.
\bs

{\bf Notation:} Let $W^{m,p}(\bar\O)$, where $\O$ is an open subset of 
$\text{\bf R}^N$, be the closure of $C^m(\bar \O)$ in the $W^{m,p}(\O)$-norm.
\bs

{\bf Theorem:} Let $Q=(0,1)^N$ then $W^{m,p}(\bar Q)=W^{m,p}(Q)$.
\bs

{\bf Proof:} Obviously $W^{m,p}(\bar Q)\subset W^{m,p}(Q)$. Hence
it is enough to prove the other inclusion.

Let $u\in W^{m,p}(Q)$. Then scale=dilate the variables with a factor 
$1-\e$. Then $u$ is transformed into $\tilde{u}$ and 
$\tilde{u}\in W^{m,p}(Q_{1/1-\e})$. Take restriction $\tilde{u}|_Q$.
Clearly $\tilde{u}|_Q\in W^{m,p}(\bar Q)$.
Let $\e$ take values $\{n^{-1}\}$, then $\{\tilde{u}_n\}$ a sequence 
which tends to $u$ in the norm of $W^{m,p}(Q)$.

End of Proof
\bs

{\bf Proof:} (Wannebo) The higher order Weierstrass approximation theorem.
\ss

This proof can regarded as good training ground for the ideas
discussed so far. It is not really meant for memorizing.
\ss

The key to get started with the proof is to change the setup and study 
another space.

We choose the space $u\in C^{Nm}(\bar Q)$ instead of $C^m(\bar Q)$.
\ss

Fix the special partial derivative $D^{\b}=\Pi_{i=1}^ND_i^m$. This way
the multiindex $\b$ is defined.
\ss

Define the set $S$ as the set of $2^{Nm}$ elements, each defined by
how it acts on the monomial $x^{\b}$. 
All transformations, which to each one--degree factor in $x^{\b}$,
say the factor $x_j$, substitute with either $x_j$ or $1-x_j$ in this case.
The set $S$ is all possible combinations of such transformations. 
We write $S=\{\sigma\}$
\ss

{\bf Identity:}
$$
1=\sum_{\sigma\in S}\sigma(x^{\b})
$$
-- Check!
\ss

Let $u\in C^{mN}(\bar Q)$.
Then by the Weierstrass approximation theorem, given $\e>0$, there
is a polynomial $Q_{\sigma}$, such that 

$$
||D^{\b}[\sigma(x^{\b})u]-Q_{\sigma}||_{L^{\infty}(Q)}< \e.
$$ 

Next step is to find a polynomial $P_{\sigma}$ such that

$$
||D^{\b}[\sigma(x^{\b})(u-P_{\sigma})]||_{L^{\infty}(Q)}< \e.
$$ 

In order to prove this observe that the equation

$$
D^{\b}[\sigma(x^{\b})P_{\sigma}]-Q_{\sigma}=0
$$
is solvable for $P_{\sigma}$ for every $Q_{\sigma}$.
The equation is linear and it is enough to solve
for $P_{\sigma}$ when $Q_{\sigma}=x^{\g}$ any $\g$.
Observe that

$$
D^{\b}[\sigma(x^{\b})x^{\g}]=const\cdot x^{\g} +\text{ lower order terms}.
$$

Hence it is possible by iteration to find $P_{\sigma}$.
Solve with higher orders monomials first subtract, iterate.
Since the equation is linear the general $Q_{\sigma}$
has a $P_{\sigma}$ solution.
\ss

We are interested in the case $|\a|\le m$ and then it follows automatically 
that $\a<\b$. This is the reason for studying
$C^{mN}(\bar Q)$ instead of $C^m(\bar Q)$.
\ss

Collect the results so far and use the Identity. Then

$$
||D^{\a}[u-\sum_{\sigma}\sigma(x^{\b})P_{\sigma}]||_{L^{\infty}(Q)}=
||D^{\a}[\sum_{\sigma}\sigma(x^{\b})(u-P_{\sigma})]||_{L^{\infty}(Q)}
=
$$
$$
=||\sum_{\sigma} D^{\a}[\sigma(x^{\b})(u-P_{\sigma})]||_{L^{\infty}(Q)}
\le
\sum_{\sigma} ||D^{\a}[\sigma(x^{\b})(u-P_{\sigma})]||_{L^{\infty}(Q)}
\le
$$

Use the more detailed Poincar\'e inequality (the Corollary)

$$
\le
\sum_{\sigma}||D^{\b}[\sigma(x^{\b})(u-P_{\sigma})]||_{L^{\infty}(Q)}
\le const.\cdot\e.
$$

This proves that $\Cal P$ is dense in $C^{mN}(\bar Q)$ 
in the $W^{m,\infty}(Q)$-norm.
\ss

In order to finish the proof, it only remains to prove that 
$C^{mN}(\bar Q)$ is dense in $C^m(\bar Q)$ in the $W^{m,\infty}(Q)$-norm.

Hence let $u\in C^m(\bar Q)$, general. The
space is complete in the $W^{m,\infty}(Q)$-norm. It is only needed
now to find a Cauchy sequence

$$
\{v_n\}_1^{\infty}\in C^{mN}(\bar Q)
$$
which converges in $W^{m,\infty}(Q)$-norm to $u$.
As simplification, assume that $Q$ has centre at the origin and 
side equal to 1.

Make the coordinate transformation 

$x\rightarrow (1-1/(n+1))x$, with $n$ positive integer.
We have the mapping $Q\rightarrow Q_n$, where $Q_n$ has side $1+1/n$. 
There is some convolution kernel with radial symmetry,
$\Psi(r)$, which has support in a ball with radius 1 and
with $\Psi\in C^{mN-m}$. 
Let 

$$
\Psi_n(r)=\l_n^N\Psi(\l_nr)
$$
-- a standard transformation.
Let $\l_n\rightarrow \infty$ in fast enough.
Now convolutions has such properties with respect to differentiation that

$$
v_n=\Psi_n*u_n|_Q\in C^{mN}(\bar Q).
$$
But $C^m(\bar Q)$ gives uniform continuity for the derivatives
up to order $m$ so it follows (Check!)
that $\{v_n\}$ is a Cauchy sequence
with limit $u$.
\ss

End of Proof.


\Refs

\ref
\no
\by [1] P.W. Jones
\paper Quasiconformal mappings and and extendability of functions
in Sobolev spaces
\jour Acta Math.
\vol 147
\yr 1981
\pages 71-88
\endref

\endRefs
\enddocument